\RequirePackage{fix-cm}
\documentclass[smallextended]{svjour3}       
\smartqed  
\usepackage{graphicx}
\usepackage{amsmath}
\usepackage{amsfonts}
\usepackage{amssymb}
\usepackage{array,longtable}
\usepackage{amscd}
\usepackage[cp1251]{inputenc}
\usepackage[english]{babel}
\usepackage[usenames]{color}
\usepackage{graphicx}
\graphicspath{{pictures/}}
\DeclareGraphicsExtensions{.pdf,.png,.jpg}
\usepackage[matrix,arrow,curve]{xy}
\begin{document}
\setcounter{tocdepth}{5}
\title{The problem of generalized $D$-stability in unbounded LMI regions and its computational aspects 
}

\titlerunning{The problem of generalized $D$-stability in unbounded LMI regions}        

\author{Olga Y. Kushel         \and
        Raffaella Pavani
}


\institute{O. Kushel \at
              Shanghai University, \\ Department of Mathematics, \\ Shangda Road 99, \\ 200444 Shanghai, China \\
              Tel.: +86(021)66133292\\
              \email{kushel@mail.ru}           
           \and
           R. Pavani \at
              Politecnico di Milano, \\ Department of Mathematics, \\ Piazza L. da Vinci 32, \\ 20133 Milano, Italia \\
              Tel.: +39(02)23994516\\
              \email{raffaella.pavani@polimi.it}
}

\date{Received: date / Accepted: date}

\maketitle

\begin{abstract}
We generalize the concepts of $D$-stability and additive $D$-stability of matrices. For this, we consider a family of unbounded regions defined in terms of Linear Matrix Inequalities (so-called LMI regions). We study the problem when the localization of a matrix spectrum in an unbounded LMI region is preserved under specific multiplicative and additive perturbations of the initial matrix. The most well-known particular cases of unbounded LMI regions (namely, conic sectors and shifted halfplanes) are considered. A new $D$-stability criterion as well as sufficient conditions for generalized $D$-stability are analyzed. Several applications of the developed theory to dynamical systems are shown.

\end{abstract}

\tableofcontents

\section{Introduction}
The concept of $D$-stability was introduced a long time ago, initially in
the papers on mathematical economics. Then in the last fifty years,
$D$-stable matrices were involved in many applications, not only in economics,
but also in control and system theory, neural networks, large scale systems,
mathematical ecology, etc. A rich list of references related to results and
applications can be found in \cite{KU2}.
Many problems involved in dynamical systems equilibria can be
treated by means of $D$-stability. About this last topic, Section 6
will present further contributions.

Firstly, we introduce three main definitions. Let $\mathcal{M}^{n\times n}$ be
the set of all square real\ $n\times n$ matrices; $\sigma (\mathbf{A})$ be the
spectrum of a matrix $\mathbf{A}$ $\in \mathcal{M}^{n\times n}$ (i.e. the set of all eigenvalues of $\mathbf A$ defined as zeroes of its characteristic polynomial $f_{\mathbf A}(\lambda):= \det(\lambda{\mathbf I - {\mathbf A}})$); ${\mathcal D}^+ \subset \mathcal{M}^{n\times n}$ be the set of all positive diagonal matrices (i.e. matrices with
positive entries on the principal diagonal while the entries outside the principal diagonal are zeroes) and ${\mathcal D}^+_0 \subset \mathcal{M}^{n\times n}$ be the set of all nonnegative diagonal matrices (i.e. diagonal matrices with
nonnegative entries on the principal diagonal).

\textbf{Definition 1} (see, for example, \cite{BELL}, \cite{KAB}). A matrix $\mathbf{A}$ $\in \mathcal{M}^{n\times n}$ is
called {\it Hurwitz stable} or just {\it stable} if ${\rm Re}(\lambda )<0$ for all $\lambda \in \sigma (\mathbf{A})$.

\textbf{Definition 2} (see \cite{AM}). A matrix $\mathbf{A}$ $\in \mathcal{M}^{n\times n}$ is
called {\it (multiplicative) $D$-stable} if ${\rm Re}(\lambda )<0$ for all $\lambda \in
\sigma (\mathbf{DA})$, where ${\mathbf D}$ is any matrix from ${\mathcal D}^+$.

\textbf{Definition 3} (see \cite{CROSS}). A matrix $\mathbf{A}$ $\in \mathcal{M}^{n\times n}$ is
called {\it (additive) $D$-stable} if ${\rm Re}(\lambda )<0$ for all $\lambda \in
\sigma (-{\mathbf D} + {\mathbf A})$, where ${\mathbf D}$ is any matrix from ${\mathcal D}_0^+$.

In practice, we may consider so-called positive stability ($\mathbf{A}$ $\in \mathcal{M}^{n\times n}$ is called {\it positively stable} if $-{\mathbf A}$ is Hurwitz stable). In this case, $D$-stability is the positive stability that is not lost under
any multiplicative (or additive) perturbation with a positive (nonnegative) diagonal matrix.

Obviously, for any stable matrix ${\mathbf A} \in {\mathcal M}^{n \times n}$ checking the property of $D$-stability via Definitions 2 and 3 is intractable, since it implies to check the negativity of the real parts of
eigenvalues for an infinite number of positive parameters appearing on the
principal diagonal of the matrix ${\mathbf D}$. Even symbolic computation of the
eigenvalues of the matrix ${\mathbf D}{\mathbf A}$ (respectively, ${\mathbf D}+{\mathbf A}$) reveals unfeasible, since for increasing
$n$ their symbolic expressions become too complicated to be
checked. Therefore, the property of $D$-stability has to be characterized by
means of theoretical necessary and sufficient conditions. Clearly, such
necessary and sufficient conditions have to be simple enough to be tractable.

The problem of $D$-stability characterization enjoys a large literature. Here we
recall that for $n\times n$ matrices with $n\leq 3$ the characterization is
relatively simple (e.g. \cite{CA}). However, for
$n\geq 4$ the characterization problem appears considerably more complicated. We
note that \cite{JOHN4} provides necessary
and sufficient conditions for the multiplicative $D$-stability of a $4\times 4$ matrix
resorting to the Routh--Hurwitz criterion, but the method appears too complicated
to be implemented efficiently within a numerical context. Indeed, later
others followed the same path (e.g. \cite{IJP},
\cite{BUR}) without succeeding in reducing the
numerical complexity of the proposed characterization.

Recently (see \cite{PAV1}, \cite{PAV2}) a new approach was proposed resorting to numerical linear algebra.
In \cite{PAV2}, the computation of a convenient
characteristic polynomial \ is used to simplify the numerical implementation
of the proposed method which allows to characterize $n\times n$
matrices with $n=2,3,4,5,..$ and more, easily. Furthermore, within the context of
LMI (Linear Matrix Inequalities) regions, introduced in \cite{CHG} and deeply studied in \cite{KU2}, a new point of view on $D$-stability is
assumed. The generalized theory of stability (with respect to LMI regions) now allows this work to present a new theory of $D$-stability
which generalizes Definition 2 (actually, the left-hand side of the complex
 plane can be also seen as a LMI region). Furthermore new and complete results are
presented. They include generalized $D$-stability criteria using additive
compound matrices also.

The paper is organized as follows. Section 2 introduces the concepts of stability and $D$-stability within the
context of LMI regions. Section 3 is devoted to the general "forbidden
boundary" approach and supplies necessary and sufficient conditions for $D$%
-stability characterization from this new point of view. Then in Section 4
additive compound matrices are used to present a new $D$-stability criterion and to generalize it. In
Section 5, the computational aspects are considered and supported by numerical
examples. At last, applications to dynamical systems are studied in details
in Section 6. Section 7 reports some concluding remarks.

\section{Stability and $D$-stability with respect to LMI regions}
\subsection{LMI regions and matrix spectra localization}

Recall the following definition (see \cite{CHG}, \cite{CGA}). A subset ${\mathfrak D} \subset {\mathbb C}$ that can be defined as
\begin{equation}\label{LMI} {\mathfrak D} = \{z \in {\mathbb C}: \ {\mathbf L} + {\mathbf M}z+{\mathbf M}^T\overline{z} \prec 0\},\end{equation}
where ${\mathbf L}, {\mathbf M} \in {\mathcal M}^{n \times n}$, ${\mathbf L}^T = {\mathbf L}$, is called an {\it LMI region} with the {\it characteristic function} $f_{\mathfrak D}(z) = {\mathbf L} + z{\mathbf M}+\overline{z}{\mathbf M}^T$ and {\it generating matrices} $\mathbf M$ and $\mathbf L$.
The above definition implies, that {\it an LMI region is open} (see \cite{KU}). In the sequel, we will use the notations $\overline{{\mathfrak D}}$ for the closure of $\mathfrak D$ and $\partial({\mathfrak D})$ for the boundary of $\mathfrak D$.

 Well-known examples of LMI regions are the left-hand side of the complex plane $${\mathbb C}^- = \{\lambda \in {\mathbb C} : {\rm Re}(\lambda) < 0\},$$
 with the characteristic function $$f_{{\mathbb C}^-}(z) = z + \overline{z};$$
 the shifted half-plane with the shift parameter $\alpha \in {\mathbb R}$
 $${\mathbb C}^-_{\alpha} = \{\lambda \in {\mathbb C} : {\rm Re}(\lambda) < \alpha\},$$
 with the characteristic function $$f_{{\mathbb C}^-_{\alpha}}(z) = z + \overline{z}-2\alpha;$$
 and the conic sector around the negative direction of the real axis with the apex at the origin and inner angle $2\theta$, $0 < \theta < \frac{\pi}{2}$,
$$
{\mathbb C}^-_{\theta} = \{z = x+iy \in {\mathbb C}: x < 0; -x\tan\theta < y < x\tan\theta\},
$$
with the characteristic function
$$f_{{\mathbb C}^-_{\theta}} = \begin{pmatrix} \sin(\theta) & \cos(\theta) \\ - \cos(\theta) & \sin(\theta) \\ \end{pmatrix}z + \begin{pmatrix} \sin(\theta) & -\cos(\theta) \\  \cos(\theta) & \sin(\theta) \\ \end{pmatrix}\overline{z}.$$

{\bf Definition 4.} Given an LMI region $\mathfrak D$, an $n \times n$ matrix ${\mathbf A}$ is called {\it stable with respect to $\mathfrak D$} or simply {\it $\mathfrak D$-stable} if $\sigma({\mathbf A}) \subset {\mathfrak D}$.

 The most well-known particular case of $\mathfrak D$-stability is Hurwitz stability (${\mathfrak D} = {\mathbb C}^-$). Here, we will also consider the case of relative (sector) stability: for a given value $\theta$, $0 < \theta < \frac{\pi}{2}$, an $n \times n$ real matrix $\mathbf A$ is called {\it relatively stable} if $\sigma({\mathbf A}) \subset {\mathbb C}^-_{\theta}$ (see, for example, \cite{DBE}, \cite{GUJU}).

  In the sequel, we will also use the following concept (see \cite{BHPM}). Given an LMI region $\mathfrak D$, an $n \times n$ matrix ${\mathbf A}$ is called {\it $\partial({\mathfrak D})$-regular} if $\sigma({\mathbf A}) \cap \partial({\mathfrak D}) = \emptyset$.

 The following generalization of the famous Lyapunov theorem is proved in \cite{CGA} (see \cite{CGA}, p. 360, Theorem 2.2, for the definition and properties of the Kronecker (or tensor) product of matrices, denoted by $\otimes$, see \cite{MA}, p. 326, also see \cite{GUJU}).
\begin{theorem}[Lyapunov theorem for LMI regions]
Given an LMI region $\mathfrak D$, defined by \eqref{LMI}, a matrix $\mathbf A$ is $\mathfrak D$-stable if and only if there is a symmetric positive definite matrix $\mathbf H$ such that the matrix
    \begin{equation}\label{LMIeq}{\mathbf W} :=\ {\mathbf L}\otimes {\mathbf H} + {\mathbf M}\otimes({\mathbf H}{\mathbf A})+{\mathbf M}^T\otimes({\mathbf A}^T{\mathbf H})\end{equation}
    is negative definite.
\end{theorem}

The particular cases of Equation \eqref{LMIeq} for ${\mathfrak D} = {\mathbb C}^-$, ${\mathfrak D} = {\mathbb C}^-_{\alpha}$ and ${\mathfrak D} = {\mathbb C}^-_{\theta}$ are as follows:
\begin{equation}\label{LyapHur} {\mathbf W} = {\mathbf H}{\mathbf A} + {\mathbf A}^T{\mathbf H}; \end{equation}
\begin{equation}{\mathbf W} = {\mathbf H}{\mathbf A} + {\mathbf A}^T{\mathbf H} - 2\alpha{\mathbf H}; \end{equation}
\begin{equation}\label{sect}{\mathbf W} =\begin{pmatrix} \sin(\theta){\mathbf H}{\mathbf A} & \cos(\theta){\mathbf H}{\mathbf A} \\ - \cos(\theta){\mathbf H}{\mathbf A} & \sin(\theta){\mathbf H}{\mathbf A} \\ \end{pmatrix} + \begin{pmatrix} \sin(\theta){\mathbf A}^T{\mathbf H} & -\cos(\theta){\mathbf A}^T{\mathbf H} \\  \cos(\theta){\mathbf A}^T{\mathbf H} & \sin(\theta){\mathbf A}^T{\mathbf H} \\ \end{pmatrix}, \end{equation}
respectively.

\subsection{Links between matrix stability and $\mathfrak D$-stability}
The Generalized Lyapunov Theorem for ${\mathfrak D} = {\mathbb C}^-_{\alpha}$ easily implies the following eigenvalue localization criterion.
\begin{theorem}\label{shift} Given a matrix ${\mathbf A} \in {\mathcal M}^{n \times n}$, and a value $\alpha \in {\mathbb R}$, the following conditions are equivalent.
\begin{enumerate}
\item[\rm (i)] $\mathbf A$ is ${\mathbb C}^-_{\alpha}$-stable;
\item[\rm (ii)] ${\mathbf A} - \alpha{\mathbf I}$ is Hurwitz stable.
\end{enumerate}
\end{theorem}

For ${\mathfrak D} = {\mathbb C}^-_{\theta}$, the following eigenvalue localization criterion was proved in \cite{DAR} (see also \cite{ABJ} for a simpler proof).

\begin{theorem}\label{tilde} Given a matrix ${\mathbf A} \in {\mathcal M}^{n \times n}$ and an angle $\theta$, $0 < \theta \leq \frac{\pi}{2}$. Then the following conditions are equivalent.
\item[\rm (i)] $\mathbf A$ is ${\mathbb C}^-_{\theta}$-stable;
\item[\rm (ii)] An $2n \times 2n$ matrix $\widetilde{{\mathbf A}}= \begin{pmatrix} \sin(\theta){\mathbf A} & -\cos(\theta){\mathbf A} \\  \cos(\theta){\mathbf A} & \sin(\theta){\mathbf A} \\ \end{pmatrix}$ is (Hurwitz) stable.
\end{theorem}

Now recall the following definitions and facts from matrix theory (see, for example, \cite{HERK2}). A matrix ${\mathbf A} \in {\mathcal M}^{n \times n}$ is called a {\it $Q$-matrix} if the inequality
$$\sum_{(i_1, \ldots, i_k)}A \begin{pmatrix}i_1 & \ldots & i_k \\ i_1 & \ldots & i_k \end{pmatrix} > 0, \qquad 1 \leq i_1 < \ldots < i_k \leq n, $$
holds for all $k, \ 1 \leq k \leq n$. The condition "$-{\mathbf A}$ is a $Q$-matrix" is necessary but not sufficient for ${\mathbb C}^-_{\theta}$-stability $(0 < \theta \leq \frac{\pi}{2})$.

 \subsection{Robust $\mathfrak D$-stability and generalization of $D$-stability}
Here, we use the notation ${\mathcal D}^+_{\geq 1}$ for the following subclass of the class of positive diagonal matrices ${\mathcal D}^+$:
$${\mathcal D}^+_{\geq 1} = \{{\mathbf D} = {\rm diag}\{d_{11}, \ \ldots, \ d_{nn}\}: 1 \leq d_{ii} < + \infty, \ i = 1, \ \ldots, \ n\}. $$

 First, we consider the generalization of Definition 2 (multiplicative $D$-stability) to the case of an arbitrary unbounded LMI region. We consider separately the cases when $0 \in \overline{{\mathfrak D}}$ and when $0 \in {\mathbb C}\setminus \overline{{\mathfrak D}}$. The following definitions make sense.

 {\bf Definition 5.} Given an unbounded LMI region $\mathfrak D$ with $0 \in \overline{{\mathfrak D}}$, we say that an $n \times n$ real matrix $\mathbf A$ is {\it (multiplicative) $D$-stable with respect to $\mathfrak D$} or simply {\it (multiplicative) $({\mathfrak D}, D)$-stable} if $\sigma({\mathbf D}{\mathbf A}) \subset {\mathfrak D}$ for every ${\mathbf D} \in {\mathcal D}^+$.

 {\bf Definition 5'.} Given an unbounded LMI region $\mathfrak D$ with $0 \in {\mathbb C} \setminus \overline{{\mathfrak D}}$, we say that an $n \times n$ real matrix $\mathbf A$ is {\it (multiplicative) $D$-stable with respect to $\mathfrak D$} or simply {\it (multiplicative) $({\mathfrak D}, D)$-stable} if $\sigma({\mathbf D}{\mathbf A}) \subset {\mathfrak D}$ for every ${\mathbf D} \in {\mathcal D}^+_{\geq 1}$.

 Now we consider the following two particular cases of LMI regions that are chosen due to a lot of applications.

 {\bf 1. Shifted half-plane ${\mathbb C}^-_{\alpha}$.} Here, we consider separately the case, when $\alpha \geq 0$ which implies $0 \in \overline{{\mathbb C}^-_{\alpha}}$ and $\alpha < 0$ which implies $0 \in {\mathbb C} \setminus \overline{{\mathbb C}^-_{\alpha}}$.

 {\bf Definition 6.} We say that an $n \times n$ real matrix $\mathbf A$ is {\it (multiplicative) $({\mathbb C}^-_{\alpha}, D)$-stable} or {\it $\alpha$-shift stable} with  $\alpha \geq 0$ if ${\rm Re}(\lambda )<\alpha$ for all $\lambda \in
\sigma (\mathbf{DA})$, where ${\mathbf D}$ is any matrix from ${\mathcal D}^+$.

{\bf Definition 6'.} We say that an $n \times n$ real matrix $\mathbf A$ is {\it (multiplicative) $({\mathbb C}^-_{\alpha}, D)$-stable} or {\it $\alpha$-shift stable} with  $\alpha < 0$ if ${\rm Re}(\lambda )<\alpha$ for all $\lambda \in
\sigma (\mathbf{DA})$, where ${\mathbf D}$ is any matrix from ${\mathcal D}^+_{\geq 1}$.

 {\bf 2. Open conic sector ${\mathbb C}^-_{\theta}$.} For this case, we have $0 \in \overline{{\mathbb C}^-_{\theta}}$ and
  obtain the following definition.

{\bf Definition 7.} For a given value $\theta$, $0 < \theta < \frac{\pi}{2}$, we call an $n \times n$ real matrix $\mathbf A$ {\it relatively $D$-stable} if $\sigma({\mathbf D}{\mathbf A}) \subset {\mathbb C}^-_{\theta}$ for every positive diagonal matrix $\mathbf D$.

Let us also introduce the following concepts, based on the notion of $\partial({\mathfrak D})$-regularity. Given an unbounded LMI region $\mathfrak D$ with $0 \in \overline{{\mathfrak D}}$ (with $0 \in {\mathbb C} \setminus \overline{{\mathfrak D}}$), we say that an $n \times n$ real matrix $\mathbf A$ is {\it (multiplicative) $(\partial{\mathfrak D}, D)$-regular} if $\sigma({\mathbf D}{\mathbf A}) \cap \partial{\mathfrak D} = \emptyset$ for every positive diagonal matrix $\mathbf D$ (respectively, for every diagonal matrix ${\mathbf D} \in {\mathcal D}^+_{\geq 1}$). The particular case of $(\partial{\mathfrak D}, D)$-regularity, based on Hurwitz stability (${\mathfrak D} = {\mathbb C}^-$, $\partial({\mathfrak D}) = OY$) is $D$-hyperbolicity: an $n \times n$ real matrix $\mathbf A$ is called {\it $D$-hyperbolic} if the eigenvalues of ${\mathbf D}{\mathbf A}$ have nonzero real parts for every real nonsingular $n \times n$ diagonal matrix $\mathbf D$ (see \cite{AB1}, also \cite{AB2}).

Next, we are interested in LMI regions with the following property.

{\bf Property 1.} For every $z \in {\mathfrak D}$, any half-line of the form $\{z - t\}_{t \geq 0}$ also lies in $\mathfrak D$, in other words, the negative direction ${\mathbb R}^-$ of the real axis $\mathbb R$ is a direction of recession of $\mathfrak D$ (see \cite{KU2}).

It was shown in \cite{KU2} (see \cite{KU2}, Theorem 11), that an LMI region $\mathfrak D$ has Property 1 if and only if its generating matrix $\mathbf M$ in Formula \eqref{LMI} is positive semidefinite (in the sense that its symmetric part ${\mathbf M} + {\mathbf M}^T$ is positive semidefinite). Thus the regions ${\mathbb C}^-_{\theta}$, $0 < \theta \leq \frac{\pi}{2}$ and ${\mathbb C}^-_{\alpha}$, $\alpha \in {\mathbb R}$ have Property 1.

Now let us consider the generalization of Definition 3 (additive $D$-stability).

{\bf Definition 8.} Given an unbounded LMI region $\mathfrak D$ with Property 1, we say that an $n \times n$ real matrix $\mathbf A$ is {\it additive $D$-stable with respect to $\mathfrak D$} or simply {\it additive $({\mathfrak D}, D)$-stable} if $\sigma({\mathbf A} - {\mathbf D}) \subset {\mathfrak D}$ for every ${\mathbf D} \in {\mathcal D}_0^+$.

The particular cases of this definition are the following.

 {\bf Definition 9.} We say that an $n \times n$ real matrix $\mathbf A$ is {\it additive $({\mathbb C}^-_{\alpha}, D)$-stable} or {\it additive $\alpha$-shift stable} if ${\rm Re}(\lambda )<\alpha$ for all $\lambda \in
\sigma ({\mathbf A} - {\mathbf D})$, where ${\mathbf D}$ is any matrix from ${\mathcal D}_0^+$.

{\bf Definition 10.} For a given value $\theta$, $0 < \theta < \frac{\pi}{2}$, we call an $n \times n$ real matrix $\mathbf A$ {\it relatively additive $D$-stable} if $\sigma({\mathbf A} - {\mathbf D}) \subset {\mathbb C}^-_{\theta}$ for every nonnegative diagonal matrix $\mathbf D$.

Comparing the obtained above definitions of $({\mathfrak D}, D)$-stability with the definitions of robust ${\mathfrak D}$-stability (see, for example, \cite{CGA}, \cite{BARMI}), we conclude that $({\mathfrak D}, D)$-stability describes a specific type of robustness with norm-unbounded uncertainties.

\subsection{Links between $D$-stability and $({\mathfrak D}, D)$-stability}

Some of the results, describing $({\mathfrak D}, D)$-stability can be easily deduced from the corresponding results for $D$-stability.

\begin{lemma}[Elementary properties of $({\mathfrak D}, D)$-stable matrices]
Given an unbounded LMI region $\mathfrak D$, let ${\mathbf A} \in {\mathcal M}^{n \times n}$ be (multiplicative) $({\mathfrak D}, D)$-stable. Then each of the following matrices is also $({\mathfrak D}, D)$-stable:
\begin{enumerate}
\item[\rm (i)] ${\mathbf A}^T$;
\item[\rm (ii)] ${\mathbf P}^T{\mathbf A}{\mathbf P}$, where $\mathbf P$ is a permutation matrix;
\item[\rm (iii)] ${\mathbf D}{\mathbf A}{\mathbf E}$, where ${\mathbf D}, \ {\mathbf E} \in {\mathcal D}^+$ if $0 \in \overline{{\mathfrak D}}$ and ${\mathbf D}, \ {\mathbf E} \in {\mathcal D}^+_{\geq 1}$ if $0 \in {\mathbb C} \setminus \overline{{\mathfrak D}}$.
\end{enumerate}
\end{lemma}

{\bf Proof.} (i) and (ii) copies the reasoning of the corresponding results of \cite{JOHN1} (see \cite{JOHN1}, p. 54, Observation (ii)) for $D$-stable matrices.

(iii) The matrix ${\mathbf D}{\mathbf A}{\mathbf E}$ is similar to the matrix ${\mathbf E}{\mathbf D}{\mathbf A}$. Obviously,  ${\mathbf E}{\mathbf D} \in {\mathcal D}^+$ if ${\mathbf E}, \ {\mathbf D} \in {\mathcal D}^+$ and ${\mathbf E}{\mathbf D} \in {\mathcal D}^+_{\geq 1}$ if ${\mathbf E}, \ {\mathbf D} \in {\mathcal D}^+_{\geq 1}$. Thus the matrix  ${\mathbf E}{\mathbf D}{\mathbf A}$ is obviously $\mathfrak D$-stable and $({\mathfrak D}, D)$-stable for every $({\mathfrak D}, D)$-stable matrix $\mathbf A$. $\square$

Note that in the case of an arbitrary unbounded LMI region $\mathfrak D$, the property of $\mathfrak D$-stability of $\mathbf A$ does not imply $\mathfrak D$-stability of ${\mathbf A}^{-1}$. However, this is true for ${\mathfrak D} = {\mathbb C}^-_{\theta}$, $0 < \theta \leq \frac{\pi}{2}$. The following result holds.

\begin{lemma} Let an invertible matrix ${\mathbf A} \in {\mathcal M}^{n \times n}$ be (multiplicative) $({\mathbb C}^-_{\theta}, D)$-stable $(0 < \theta \leq \frac{\pi}{2})$. Then ${\mathbf A}^{-1}$ is also $({\mathbb C}^-_{\theta}, D)$-stable.
\end{lemma}

The proof copies the reasoning of \cite{JOHN1} (see \cite{JOHN1}, p. 54, Observation (ii)).

Now let us recall the following definitions (see \cite{HERK2}, \cite{CROSS}). A matrix ${\mathbf A} \in {\mathcal M}^{n \times n}$ is called a {\it $P$-matrix ($P_0$-matrix)} if all its principal minors are positive (respectively, nonnegative), i.e the inequality $A \left(\begin{array}{ccc}i_1 & \ldots & i_k \\ i_1 & \ldots & i_k \end{array}\right) > 0$ (respectively, $\geq 0$)
holds for all $(i_1, \ \ldots, \ i_k), \ 1 \leq i_1 < \ldots < i_k \leq n$, and all $k, \ 1 \leq k \leq n$. A matrix ${\mathbf A} \in {\mathcal M}^{n \times n}$ is called a {\it $P_0^+$-matrix} if it is a $P_0$-matrix and, in addition, the sums of all principal minors of every fixed order $i$ are positive $(i = 1, \ \ldots, \ n)$ (i.e. if it is a $P_0$-matrix and a $Q$-matrix at the same time).

The following necessary condition for $D$-stability was established in \cite{CROSS} (see \cite{CROSS}, p. 256, Corollary): {\it if $\mathbf A$ is (multiplicative) $D$-stable then $\mathbf A$ is a $P_0^+$-matrix.}

\section{"Forbidden boundary" approach}
\subsection{General "forbidden boundary" approach}
Here, we state the following general result.
\begin{theorem}[Forbidden boundary]\label{fobou}
Let a matrix ${\mathbf A} \in {\mathcal M}^{n \times n}$ be $\mathfrak D$-stable with respect to an unbounded LMI region $\mathfrak D$. Then the following conditions are equivalent.
\begin{enumerate}
\item[\rm (i)] $\mathbf A$ is (multiplicative) $({\mathfrak D}, D)$-stable;
\item[\rm (ii)] $\mathbf A$ is ${\mathfrak D}$-stable and $(\partial({\mathfrak D}), D)$-regular.
\end{enumerate}
\end{theorem}
{\bf Proof.} $\Rightarrow$ Since ${\mathbf I} \in {\mathcal D}^+_{\geq 1} \subset {\mathcal D}^+$, we conclude, that $({\mathfrak D}, D)$-stability implies $\mathfrak D$-stability. Since LMI regions are open, $({\mathfrak D}, D)$-stability implies also $(\partial({\mathfrak D}), D)$-regularity.

$\Leftarrow$ Conversely, let $\mathbf A$ be ${\mathfrak D}$-stable and $(\partial({\mathfrak D}), D)$-regular. From this, we conclude that ${\mathbf D}{\mathbf A}$ does not have any eigenvalues on $\partial({\mathfrak D})$ whenever ${\mathbf D}$ is positive diagonal (respectively, whenever ${\mathbf D} \in {\mathcal D}^+_{\geq 1}$). Assume there is ${\mathbf D}_0 \in {\mathcal D}^+$ (respectively, ${\mathbf D}_0 \in {\mathcal D}^+_{\geq 1}$) such that ${\mathbf D}_0{\mathbf A}$ is not $\mathfrak D$-stable. Then for the parameter-dependent matrix family $\{{\mathbf D}(t)\}$, ${\mathbf D}(t) = t{\mathbf D}_0 + (1-t){\mathbf I}$, $t \in [0,1]$, we get ${\mathbf D}(t) \in {\mathcal D}^+$ (respectively, ${\mathbf D}(t) \in {\mathcal D}^+_{\geq 1}$) for every $t \in [0,1]$. Then, ${\mathbf D}(0){\mathbf A}= {\mathbf A}$ is $\mathfrak D$-stable and ${\mathbf D}(1){\mathbf A} = {\mathbf D}_0{\mathbf A}$ is not $\mathfrak D$-stable. Since the eigenvalues of ${\mathbf D}(t){\mathbf A}$ are changing continuously on $t$, we have $t_0 \in [0,1]$ such that ${\mathbf D}(t_0){\mathbf A}$ have eigenvalues on the boundary $\partial({\mathfrak D})$ of $\mathfrak D$. This contradicts $(\partial({\mathfrak D}), D)$-regularity of $\mathbf A$.
 $\square$

The theorem will remain valid if we replace ${\mathfrak D}$-stability with $({\mathbb C}\setminus\overline{{\mathfrak D}})$-stability.

\subsection{"Forbidden boundary" conditions for conic regions}
Here, we consider the regions bounded by the rays defined by $y = x \tan\theta$ and $y = -x \tan \theta$, $x \in [0, \infty)$, i.e. the conic sector ${\mathbb C}^-_{\theta}$, $0  < \theta \leq \frac{\pi}{2}$, and the interior part of its completion ${\mathbb C} \setminus \overline{{\mathbb C}^-_{\theta}}$.

\begin{theorem}\label{bouncon} Let ${\mathbf A} \in {\mathcal M}^{n \times n}$ be nonsingular, and $z = \cos\theta + i\sin\theta$ be a fixed point on the unit circle. The following statements are equivalent: "for all positive diagonal matrices $\mathbf D$, ..."
\begin{enumerate}
\item[\rm (i)] ${\mathbf D}{\mathbf A}$ has no eigenvalues on the rays defined by $\{tz\}_{t = 0}^{\infty}$ and $\{t\overline{z}\}_{t = 0}^{\infty}$.
\item[\rm (ii)] Both ${\mathbf A} - z{\mathbf D}$ and ${\mathbf A} - \overline{z}{\mathbf D}$ are nonsingular.
\item[\rm (iii)] $({\mathbf D}{\mathbf A})^2 - 2\cos\theta{\mathbf D}{\mathbf A}$ does not have any eigenvalues equal to $- 1$ .
\item[\rm (iv)] ${\mathbf A}{\mathbf D}^{-1} + {\mathbf D}{\mathbf A}^{-1} - 2\cos\theta{\mathbf I}$ is nonsingular.
\item[\rm (v)] $\det \begin{pmatrix}{\mathbf A} & &{\mathbf D} \\ -{\mathbf D} & & {\mathbf A} - 2\cos\theta{\mathbf D} \end{pmatrix} \neq 0$.
\item[\rm (vi)] The characteristic polynomial $f_{{\mathbf D}{\mathbf A}}(\lambda)$ of ${\mathbf D}{\mathbf A}$ is not divisible by the polynomial $\lambda^2 - 2\cos\theta \lambda + 1$.
\end{enumerate}
\end{theorem}
{\bf Proof.} ${\rm (i)} \Rightarrow {\rm (ii)}$. Suppose that $\det({\mathbf A} - z{\mathbf D}) = 0$ for some positive diagonal matrix $\mathbf D$. Since $\det({\mathbf A} - z{\mathbf D}) = \det({\mathbf D})\det({\mathbf D}^{-1}{\mathbf A} - z{\mathbf I})$, we have that $z \in \sigma({\mathbf D}^{-1}{\mathbf A})$. Thus, for some positive diagonal matrix $\widetilde{{\mathbf D}}:={\mathbf D}^{-1}$, the matrix $\widetilde{{\mathbf D}}{\mathbf A}$ has an eigenvalue $z \in \{tz\}_{t = 0}^{\infty}$. The inequality $\det({\mathbf A} - \overline{z}{\mathbf D}) \neq 0$ is proved analogically.

${\rm (ii)} \Rightarrow {\rm (iii)}$. Suppose that $-1 \in \sigma(({\mathbf D}{\mathbf A})^2 - 2\cos\theta{\mathbf D}{\mathbf A})$ for some positive diagonal matrix $\mathbf D$. Then $\det(({\mathbf D}{\mathbf A})^2 - 2\cos\theta{\mathbf D}{\mathbf A} + {\mathbf I}) = 0$ together with $z = \cos\theta + i\sin\theta$ imply $$\det(({\mathbf D}{\mathbf A} - z{\mathbf I})({\mathbf D}{\mathbf A} - \overline{z}{\mathbf I})) = \det(({\mathbf D}{\mathbf A} - z{\mathbf I})\det({\mathbf D}{\mathbf A} - \overline{z}{\mathbf I}))=0,$$
which contradicts ${\rm (ii)}$.

${\rm (iii)} \Rightarrow {\rm (iv)}$. Obviously, $$\det({\mathbf A}{\mathbf D}^{-1} + {\mathbf D}{\mathbf A}^{-1}- 2\cos\theta{\mathbf I}) = \det({\mathbf D}{\mathbf A}^{-1}) \det(({\mathbf A}{\mathbf D}^{-1})^2 + {\mathbf I} - 2\cos\theta{\mathbf A}{\mathbf D}^{-1}).$$
From ${\rm (iii)}$ and the similarity transformation ${\mathbf D}^{-1}(({\mathbf A}{\mathbf D}^{-1})^2 - 2\cos\theta{\mathbf A}{\mathbf D}^{-1} ){\mathbf D} = ({\mathbf D}^{-1}{\mathbf A})^2 - 2\cos\theta{\mathbf D}^{-1}{\mathbf A}$, we get $\det(({\mathbf A}{\mathbf D}^{-1})^2 + {\mathbf I}  - 2\cos\theta{\mathbf A}{\mathbf D}^{-1}) \neq 0$.
 By observing that $\det({\mathbf A}^{-1}{\mathbf D}) = \frac{\det({\mathbf D})}{\det({\mathbf A})} \neq 0$, we complete the proof.

${\rm (iv)} \Rightarrow {\rm (v)}$. By Schur completion formula (see, e.g. \cite{BER}, p. 135),
$$\det \begin{pmatrix}{\mathbf A} & &{\mathbf D} \\ -{\mathbf D} & &{\mathbf A} - 2\cos\theta{\mathbf D} \end{pmatrix} = \det({\mathbf A})\det({\mathbf A} - 2\cos\theta{\mathbf D} + {\mathbf D}{\mathbf A}^{-1}{\mathbf D}) = $$ $$ \det({\mathbf A})\det({\mathbf D})\det({\mathbf A}{\mathbf D}^{-1} + {\mathbf D}{\mathbf A}^{-1} - 2\cos\theta{\mathbf I}) \neq 0.$$

${\rm (v)} \Rightarrow {\rm (vi)}$. Continuing the above equality, we get $$\det \begin{pmatrix}{\mathbf A} & {\mathbf D} \\ -{\mathbf D} & {\mathbf A}- 2\cos\theta{\mathbf D}  \end{pmatrix} = $$ $$\det({\mathbf A})\det({\mathbf D})\det({\mathbf D}{\mathbf A}^{-1})\det(({\mathbf A}{\mathbf D}^{-1})^2 -2\cos\theta {\mathbf A}{\mathbf D}^{-1} + {\mathbf I}) \neq 0,$$
which implies $\det(({\mathbf A}{\mathbf D}^{-1})^2 -2\cos\theta {\mathbf A}{\mathbf D}^{-1} + {\mathbf I}) \neq 0$. Re-writing $\det(({\mathbf A}{\mathbf D}^{-1})^2 -2\cos\theta {\mathbf A}{\mathbf D}^{-1} + {\mathbf I})$ in the form of the product of eigenvalues of the corresponding matrix, we get $$\det(({\mathbf A}{\mathbf D}^{-1})^2 -2\cos\theta {\mathbf A}{\mathbf D}^{-1} + {\mathbf I}) =\prod_{i=1}^n(1 -2\cos\theta\lambda_i  +\lambda_i^2),$$ where $\lambda_i$, $i = 1, \ \ldots, \ n$, are the eigenvalues of the matrix ${\mathbf A}{\mathbf D}^{-1}$ (or, equivalently, $\widetilde{{\mathbf D}}{\mathbf A}$). The inequality $\prod_{i=1}^n(1-2\cos\theta\lambda_i+\lambda_i^2) \neq 0$ shows that $(1-2\cos\theta\lambda_i+\lambda_i^2) \neq 0$ whenever $\lambda_i \in \sigma(\widetilde{{\mathbf D}}{\mathbf A})$, thus $\widetilde{{\mathbf D}}{\mathbf A}$ does not have any eigenvalues equal to $z$ or $\overline{z}$ and the characteristic polynomial $f_{\widetilde{{\mathbf D}}{\mathbf A}}(\lambda)$ is not divisible by $\lambda^2-2\cos\theta\lambda + 1$.

${\rm (vi)} \Rightarrow {\rm (i)}$. Let the characteristic polynomial of ${\mathbf D}{\mathbf A}$ be not divisible by $\lambda^2-2\cos\theta\lambda + 1$ for any positive diagonal matrix $\mathbf D$. Thus ${\mathbf D}{\mathbf A}$ has no eigenvalues equal to $z$ or $\overline{z}$ for any positive diagonal matrix $\mathbf D$. Then, ${\mathbf D}{\mathbf A}$ has no eigenvalues of the form $tz$ or $t\overline{z}$, where $t > 0$, otherwise the matrix $\widehat{{\mathbf D}}{\mathbf A} = (\frac{1}{t}{\mathbf D})({\mathbf A})$ will have eigenvalues equal to $z$ or $\overline{z}$. $\square$

Combining Theorem \ref{fobou} with Theorem \ref{bouncon}, we obtain the following criterion of relative $D$-stability.

\begin{theorem} Let ${\mathbf A} \in {\mathcal M}^{n \times n}$ be ${\mathbb C}^-_{\theta}$-stable (${\mathbb C}\setminus \overline{{\mathbb C}^-_{\theta}}$-stable). Then the following conditions are equivalent.
\begin{enumerate}
\item[\rm (i)] ${\mathbf A}$ is relatively $D$-stable (respectively, $({\mathfrak D}, D)$-stable for ${\mathfrak D} = {\mathbb C}\setminus \overline{{\mathbb C}^-_{\theta}}$).
\item[\rm (ii)]  $\det({\mathbf A} - z{\mathbf D}) \neq 0$ and $\det({\mathbf A} - \overline{z}{\mathbf D}) \neq 0$ for every positive diagonal matrix $\mathbf D$ and $z = \cos\theta + i\sin\theta$.
\item[\rm (iii)] $\det({\mathbf D}{\mathbf A})^2 - 2\cos\theta{\mathbf D}{\mathbf A} + {\mathbf I}) \neq 0$ for every positive diagonal matrix $\mathbf D$.
\item[\rm (iv)] $\det({\mathbf A}{\mathbf D}^{-1} + {\mathbf D}{\mathbf A}^{-1} - 2\cos\theta{\mathbf I}) \neq 0$ for every positive diagonal matrix $\mathbf D$.
\item[\rm (v)] $\det \begin{pmatrix}{\mathbf A} & &{\mathbf D} \\ -{\mathbf D} & & {\mathbf A} - 2\cos\theta{\mathbf D} \end{pmatrix} \neq 0$ for every positive diagonal matrix $\mathbf D$.
\item[\rm (vi)] The characteristic polynomial $f_{{\mathbf D}{\mathbf A}}(\lambda)$ of ${\mathbf D}{\mathbf A}$ is not divisible by $\lambda^2 - 2\cos\theta \lambda + 1$ for every positive diagonal matrix $\mathbf D$.
\end{enumerate}
\end{theorem}

\subsection{"Forbidden boundary" conditions for the classical case ${\mathfrak D} = {\mathbb C}^-$}
The particular case of Theorem \ref{bouncon} for $\theta = \frac{\pi}{2}$, $z = \pm i$ provides the following "forbidden boundary" conditions.

\begin{theorem}\label{bounconi} Let ${\mathbf A} \in {\mathcal M}^{n \times n}$ be nonsingular. The following statements are equivalent: "for all positive diagonal matrices $\mathbf D$, ..."
\begin{enumerate}
\item[\rm (i)] ${\mathbf D}{\mathbf A}$ has no eigenvalues on the imaginary axis.
\item[\rm (ii)] ${\mathbf A} \pm i{\mathbf D}$ is nonsingular.
\item[\rm (iii)] $({\mathbf D}{\mathbf A})^2$ does not have any eigenvalues equal to $- 1$ .
\item[\rm (iv)] ${\mathbf A}{\mathbf D}^{-1} + {\mathbf D}{\mathbf A}^{-1}$ is nonsingular.
\item[\rm (v)] $\det \begin{pmatrix}{\mathbf A} & {\mathbf D} \\ -{\mathbf D} & {\mathbf A} \end{pmatrix} \neq 0$.
\item[\rm (vi)] The characteristic polynomial $f_{{\mathbf D}{\mathbf A}}(\lambda)$ of ${\mathbf D}{\mathbf A}$ is not divisible by $\lambda^2 + 1$.
\end{enumerate}
\end{theorem}

Theorem \ref{fobou} together with the equivalent conditions listed in Theorem \ref{bounconi} imply the following criteria of $D$-stability proved by different authors in \cite{JOHN5}, \cite{JOHNT}, \cite{PAV1}-\cite{PAV2}.

\begin{theorem}\label{stabcond} Let ${\mathbf A} \in {\mathcal M}^{n \times n}$ be stable. Then the following conditions are equivalent.
\begin{enumerate}
\item[\rm (i)] ${\mathbf A}$ is $D$-stable.
\item[\rm (ii)] $\det({\mathbf A} \pm i{\mathbf D}) \neq 0$ for every positive diagonal matrix $\mathbf D$.
\item[\rm (iii)] $\det(({\mathbf D}{\mathbf A})^2 + {\mathbf I}) \neq 0$ for every positive diagonal matrix $\mathbf D$.
\item[\rm (iv)] $\det({\mathbf A}{\mathbf D}^{-1} + {\mathbf D}{\mathbf A}^{-1}) \neq 0$ for every positive diagonal matrix $\mathbf D$.
\item[\rm (v)] $\det \begin{pmatrix}{\mathbf A} & {\mathbf D} \\ -{\mathbf D} & {\mathbf A} \end{pmatrix} \neq 0$ for every positive diagonal matrix $\mathbf D$.
\item[\rm (vi)] The characteristic polynomial $f_{{\mathbf D}{\mathbf A}}(\lambda)$ of ${\mathbf D}{\mathbf A}$ is not divisible by $\lambda^2 + 1$ for every positive diagonal matrix $\mathbf D$.
\end{enumerate}
\end{theorem}

 The equivalence $(i) \Leftrightarrow (ii)$ is proved in \cite{JOHN5}, p. 89, Corollary 2. The equivalence $(i) \Leftrightarrow (iv)$ is proved in \cite{PAV1}, p. 310, Corollary, using the results from \cite{JOHNT}. The equivalence $(i) \Leftrightarrow (v)$ is established by geometric methods in \cite{JOHNT}, p. 302, Proposition 2.3. The equivalence $(i) \Leftrightarrow (vi)$ is proved in \cite{PAV2}, p. 3, Proposition 2.

\section{Generalized $D$-stability criteria using additive compound matrices}
\subsection{Stability and $\mathfrak D$-stability criteria using additive compound matrices}
First, recall the following definition. Given an $n \times n$ matrix ${\mathbf A} = \{a_{ij}\}_{i,j = 1}^n$ and an $n \times n$ identity matrix ${\mathbf I} = \{\delta_{ij}\}_{i,j = 1}^n$, {\it the second additive compound matrix} ${\mathbf A}^{[2]} = \{a^{[2]}_{\alpha\beta}\}_{\alpha, \beta = (1,2)}^{(n-1, n)}$ is a matrix that consists of the sums of minors of the following form:
$$a^{[2]}_{\alpha\beta} = \begin{vmatrix} a_{ik} & \delta_{il} \\ a_{jk} & \delta_{jl} \end{vmatrix} + \begin{vmatrix} \delta_{ik} & a_{il} \\ \delta_{jk} & a_{jl} \end{vmatrix}, $$
where $\alpha = (i,j), \  1 \leq i < j \leq n$, $\beta = (k,l), \  1 \leq k < l \leq n$,
listed in the lexicographic order. The matrix  ${\mathbf A}^{[2]}$ is $\binom{n}2 \times \binom{n}2$ dimensional (see \cite{FID}).

This definition was given in \cite{FID}, but in fact, additive compound matrices (up to the inverse re-numeration of the coordinates), as well as more general constructions of exterior produts of matrices, were introduced earlier by Stephanos (see \cite{STE}) and were referred in \cite{FUL}, \cite{JUR} as "bialternate sum of ${\mathbf A}$ with itself".

A criterion of matrix stability in terms of the characteristic polynomial coefficients of the matrix and its second additive compound, was established in \cite{FUL} (see \cite{FUL}, p. 90. Theorem 9, also \cite{JUR}, p. 107, Theorem 3.15). Re-writing this criterion in terms of matrix properties, we obtain the following statement.

\begin{theorem}\label{addcomp} Let ${\mathbf A} \in {\mathcal M}^{n \times n}$, $n > 1$, and ${\mathbf A}^{[2]}$ be its second additive compound matrix. Then for $\mathbf A$ to be (Hurwitz) stable, it is necessary and sufficient that both $-{\mathbf A}$ and $-{\mathbf A}^{[2]}$ are $Q$-matrices.
\end{theorem}

Another stability criterion in terms of additive compound matrices is provided in \cite{LIW1}.

\subsection{$D$-stability criterion using additive compound matrices}
The following $D$-stability criterion can be easily deduced from Theorem \ref{addcomp}.
\begin{theorem} Let ${\mathbf A} \in {\mathcal M}^{n \times n}$, $n > 1$, and ${\mathbf A}^{[2]}$ be its second additive compound matrix. Then for $\mathbf A$ to be (multiplicative) $D$-stable, it is necessary and sufficient that the following two conditions hold:
\begin{enumerate}
\item[\rm (i)] $-{\mathbf A}$ is a $P_0^+$-matrix;
\item[\rm (ii)] $-({\mathbf D}{\mathbf A})^{[2]}$ is a $Q$-matrix for every positive diagonal matrix $\mathbf D$.
\end{enumerate}
\end{theorem}
{\bf Proof.} $\Rightarrow$ Let $\mathbf A$ be $D$-stable. Then, by the results of Cross (see \cite{CROSS}, p. 256, Corollary), $\mathbf A$ is a $P_0^+$-matrix. By the definition of $D$-stability, ${\mathbf D}{\mathbf A}$ is stable for every positive diagonal matrix $\mathbf D$. Applying Theorem \ref{addcomp} to ${\mathbf D}{\mathbf A}$, we get that $-({\mathbf D}{\mathbf A}^{[2]})$ is a $Q$-matrix for every positive diagonal matrix $\mathbf D$.

$\Leftarrow$ Let $-\mathbf A$ be a $P_0^+$-matrix. The identity $$({\mathbf D}{\mathbf A})\left(\begin{array}{ccc}i_1 & \ldots & i_k \\ i_1 & \ldots & i_k \end{array}\right) ={\mathbf D}\left(\begin{array}{ccc}i_1 & \ldots & i_k \\ i_1 & \ldots & i_k \end{array}\right){\mathbf A}\left(\begin{array}{ccc}i_1 & \ldots & i_k \\ i_1 & \ldots & i_k \end{array}\right),$$
which holds for all $(i_1, \ \ldots, \ i_k), \ 1 \leq i_1 < \ldots < i_k \leq n$, and all $k, \ 1 \leq k \leq n$, implies that $-{\mathbf D}{\mathbf A}$ is a $P_0^+$-matrix and thus a $Q$-matrix for every positive diagonal matrix $\mathbf D$. Since $-({\mathbf D}{\mathbf A})^{[2]}$ is also a $Q$-matrix, we get by Theorem \ref{addcomp} that ${\mathbf D}{\mathbf A}$ is stable for every positive diagonal matrix $\mathbf D$. This means $\mathbf A$ is $D$-stable. $\square$

Similar criterion can be easily proven for additive $D$-stability.

\subsection{$({\mathfrak D}, D)$-stability criteria using additive compound matrices}

For relative $D$-stability, we obtain the following criterion.

\begin{theorem} Given a matrix ${\mathbf A} \in {\mathcal M}^{n \times n}$, and an angle $\theta$, $0 < \theta \leq \frac{\pi}{2}$. Let $$\widetilde{{\mathbf A}}= \begin{pmatrix} \sin(\theta){\mathbf A} & -\cos(\theta){\mathbf A} \\  \cos(\theta){\mathbf A} & \sin(\theta){\mathbf A} \\ \end{pmatrix}.$$ Then for $\mathbf A$ to be (multiplicative) $({\mathbb C}^-_{\theta}, D)$-stable, it is necessary and sufficient that both $-(\widetilde{{\mathbf D}}\widetilde{{\mathbf A}})$ and $-(\widetilde{{\mathbf D}}\widetilde{{\mathbf A}})^{[2]}$ are $Q$-matrices for every block diagonal $2n \times 2n$ matrix  $$\widetilde{{\mathbf D}}= \begin{pmatrix} {\mathbf D} & {\mathbf O} \\  {\mathbf O} & {\mathbf D} \\ \end{pmatrix},$$ where $\mathbf D$ is an arbitrary  positive diagonal $n \times n$ matrix, $\mathbf O$ is a zero matrix.
\end{theorem}
{\bf Proof.} $\Rightarrow$ Let $\mathbf A$ be $({\mathbb C}^-_{\theta}, D)$-stable. Then, by the definition of $({\mathbb C}^-_{\theta}, D)$-stability, ${\mathbf D}{\mathbf A}$ is ${\mathbb C}^-_{\theta}$-stable for every positive diagonal matrix $\mathbf D$. Applying Theorem \ref{tilde} to ${\mathbf D}{\mathbf A}$, we obtain $\begin{pmatrix} \sin(\theta){\mathbf D}{\mathbf A} & -\cos(\theta){\mathbf D}{\mathbf A} \\  \cos(\theta){\mathbf D}{\mathbf A} & \sin(\theta){\mathbf D}{\mathbf A} \\ \end{pmatrix} = \widetilde{{\mathbf D}}\widetilde{{\mathbf A}}$ is (Hurwitz) stable. Then applying Theorem \ref{addcomp} to $\widetilde{{\mathbf D}}\widetilde{{\mathbf A}}$, we get that both $-(\widetilde{{\mathbf D}}\widetilde{{\mathbf A}})$ and $-(\widetilde{{\mathbf D}}\widetilde{{\mathbf A}})^{[2]}$ are $Q$-matrices for every block diagonal matrix $\widetilde{\mathbf D}$.

$\Leftarrow$ Let both $-(\widetilde{{\mathbf D}}\widetilde{{\mathbf A}})$ and $-(\widetilde{{\mathbf D}}\widetilde{{\mathbf A}})^{[2]}$ be $Q$-matrices for every block diagonal matrix $\widetilde{\mathbf D}$. Then by Theorem \ref{addcomp}, the matrix $\widetilde{{\mathbf D}}\widetilde{{\mathbf A}}$ is Hurwitz stable and by Theorem \ref{tilde}, ${\mathbf D}{\mathbf A}$ is ${\mathbb C}^-_{\theta}$-stable for every positive diagonal matrix $\mathbf D$. This condition exactly means $({\mathbb C}^-_{\theta}, D)$-stability of $\mathbf A$. $\square$

For $\alpha$-shift stability, we obtain the following criterion.

\begin{theorem} Given a matrix ${\mathbf A} \in {\mathcal M}^{n \times n}$, and a value $\alpha \in {\mathbb R}$. Then for $\mathbf A$ to be (multiplicative) $({\mathbb C}^-_{\alpha}, D)$-stable, it is necessary and sufficient that both $-({\mathbf D}{\mathbf A} - \alpha{\mathbf I})$ and $-({\mathbf D}{\mathbf A} - \alpha{\mathbf I})^{[2]}$ are $Q$-matrices for every ${\mathbf D} \in {\mathcal D}^+$ if $\alpha \geq 0$ and for every ${\mathbf D} \in {\mathcal D}^+_{\geq 1}$ if $\alpha < 0$.
\end{theorem}
{\bf Proof.} $\Rightarrow$ Let $\mathbf A$ be $({\mathbb C}^-_{\alpha}, D)$-stable. Then, by the definition of $({\mathbb C}^-_{\alpha}, D)$-stability, ${\mathbf D}{\mathbf A}$ is ${\mathbb C}^-_{\alpha}$-stable for every ${\mathbf D} \in {\mathcal D}^+$ if $\alpha \geq 0$ and for every ${\mathbf D} \in {\mathcal D}^+_{\geq 1}$ if $\alpha < 0$. Applying Theorem \ref{shift} to ${\mathbf D}{\mathbf A}$, we obtain that ${\mathbf D}{\mathbf A} - \alpha{\mathbf I}$ is (Hurwitz) stable. Then applying Theorem \ref{addcomp} to ${\mathbf D}{\mathbf A} - \alpha{\mathbf I}$, we get that both $-({\mathbf D}{\mathbf A} - \alpha{\mathbf I})$ and $-({\mathbf D}{\mathbf A} - \alpha{\mathbf I})^{[2]}$ are $Q$-matrices for every ${\mathbf D} \in {\mathcal D}^+$ if $\alpha \geq 0$ and for every ${\mathbf D} \in {\mathcal D}^+_{\geq 1}$ if $\alpha < 0$.

$\Leftarrow$ Let both $-({\mathbf D}{\mathbf A} - \alpha{\mathbf I})$ and $-({\mathbf D}{\mathbf A} - \alpha{\mathbf I})^{[2]}$ be $Q$-matrices for every ${\mathbf D} \in {\mathcal D}^+$ if $\alpha \geq 0$ and for every ${\mathbf D} \in {\mathcal D}^+_{\geq 1}$ if $\alpha < 0$. Then by Theorem \ref{addcomp}, the matrix ${\mathbf D}{\mathbf A} - \alpha{\mathbf I}$ is Hurwitz stable and by Theorem \ref{shift}, ${\mathbf D}{\mathbf A}$ is ${\mathbb C}^-_{\alpha}$-stable for every ${\mathbf D} \in {\mathcal D}^+$ if $\alpha \geq 0$ and for every ${\mathbf D} \in {\mathcal D}^+_{\geq 1}$ if $\alpha < 0$. This condition exactly means $({\mathbb C}^-_{\alpha}, D)$-stability of $\mathbf A$. $\square$

\section{Computational aspects and numerical examples}
We may think that, once we know something is \textit{computable, }whether it
takes $10$ seconds or $20$ seconds to be computed is \ a concern of
engineers. Indeed, at the present,  the
fastest computer in the world achieves the overall performance of $93.01$ \
petaflops ($10^{15}$ flops/s) so that any simple computation can be done in
an imperceptible fraction of a second. But that conclusion would not be so
obvious, if the question were \ one of $10\sec .$ versus $10^{10^{10}}\sec .$
This means that quantitative gaps have to be considered qualitative gaps
also. \ Examples are provided by the difference between reading a $400$-page
book and reading \textit{every possible such book, }or between writing down
a thousand-digit number and counting to that number.
\\
From a theoretical point of view, an algorithm is called {\it efficient}
if its running
time on inputs of size $n$ can be upper-bounded by any polynomial function
of $n$, i.e. its running time is $O(n^k)$ for some $k > 0$; an algorithm is called
{\it inefficient} if its running time can be upper bounded by $2^{poly(n)}$, where $poly(n)$ is some
polynomial in $n$; this happens, for example, when the running time is $O(2^{n^k}
)$; analogously, when
the running time is $O(n!)$ the algorithm is said to be inefficient, as well.
This choice is motivated by
some main reasons. First of all, obviously, there is a great difference between running times
of polynomial and exponential algorithms.  Then, polynomial-time algorithms are closed
under superposition. This means that if an algorithm makes polynomial
number of calls to a function that is implemented as a polynomial algorithm,
the resulting algorithm has also polynomial time-complexity. This greatly
simplifies theoretical analysis of algorithms. \smallskip

Even though the polynomial/exponential distinction is open to some objections,
empirically it reveals a useful practical tool.

About $D$-stability, we recall that in \cite{CFY} Chen et al.  established the equivalence between the $D$-stability of a real matrix and some condition on the real structured singular value of a complex matrix with a somewhat special structure. Since it was previously proved that the exact computation of a real structured singular value is in general a NP-hard problem, Chen et al. conclude that "the problem of checking $D$-stability may be generally intractable when $n$ is large".

Indeed, in \cite{JOHNT} Johnson and Tesi proposed a method which resorts to Hurwitz criterion and requires computation of Hurwitz determinants, but the computational complexity of determinant calculation is $n!$ (which means that the running time is O(n!)). Moreover, many other
authors tried to find a general algorithm to characterize $D$-stability with the
property of being simple enough to be computationally tractable, but without
any success. Instead, the approach presented here appears new and promising
in reducing the computational complexity of the $D$-stability characterization
problem. This approach is mainly based on the symbolic computation of matrices
${\mathbf D}{\mathbf A}$. If $n$ is the degree of the characteristic polynomial, any symbolic software
(such as for example MuPad) can accomplish its computation by a Hessemberg
algorithm which requires $O(n^3)$ operations, which means polynomial time.
Then, for example, for (iv) in Theorem 8, that is enough, since the determinant
is different from zero when all the eigenvalues are different from zero; this result
can be easily derived when the last coefficient of the characteristic polynomial
(which is the product of all the eigenvalues) is positive. Alternatively, (vi) in
Theorem 8 can be used, but also in this case the overall computational complexity
remains polynomial, by superposition (see \cite{Pavani 2018} for more details).

The main gain is that we transformed the problem of determinants
evaluation (which exhibits a computational complexity of $n!$) in a problem
of evaluation of some characteristic polynomial roots. Therefore, the
\textbf{global complexity appears polynomial,} and this supports the
previous conclusion that $D$-stability characterization problem \textbf{%
reveals in practice numerically tractable}.

As a numerical example we propose  the example already presented by L.A.
Burlakova \cite{BUR}, so that comparisons can be discussed. Consider the general linear
mechanical  system
\begin{equation}\label{systbur} \ddot{x} - {\mathbf B}\dot{x}-{\mathbf C}x=0, \end{equation}
where $x\in {\mathbb R}^{n}$, ${\mathbf B}, \ {\mathbf C} \in {\mathcal M}^{n \times n}$. Let ${\mathbf A} \in  {\mathcal M}^{2n\times 2n}$ be the matrix of differential equation \eqref{systbur}:

$${\mathbf A} =
\begin{pmatrix}
{\mathbf B} & {\mathbf C} \\
{\mathbf I} & {\mathbf O}
\end{pmatrix},$$
where ${\mathbf I}$ is an $n \times n$ identity matrix and $\mathbf O$ is an $n \times n$ zero matrix. Details about this
notation as well as about the relations between the property of $D$-stability of the matrix $\mathbf A$ and the properties of the related mechanical system will be provided in the next Section. In order to check the $D$-stability
of ${\mathbf A}$, Burlakova resorts to the verification of the positiveness
everywhere in the positive orthant of the Hurwitz determinants for matrix ${\mathbf D}{\mathbf A}$.  Instead,  here we propose our method. The matrix $\mathbf A$ is taken

$$ {\mathbf A} =
\begin{pmatrix}
-1 & -1 & -1 & -4/5 \\
-4 & -5 & -4 & -4 \\
1 & 0 & 0 & 0 \\
0 & 1 & 0 & 0
\end{pmatrix}, $$
as suggested in Example 1 by Burlakova. The matrix ${\mathbf A}$ is stable since the
real parts of all its eigenvalues are negative. Then, assuming $\alpha
,\beta ,\gamma ,\delta >0,$ we introduce

$${\mathbf D} = \begin{pmatrix}
\alpha  & 0 & 0 & 0 \\
0 & \beta  & 0 & 0 \\
0 & 0 & \gamma  & 0 \\
0 & 0 & 0 & \delta
\end{pmatrix}.$$

The characteristic polynomial of ${\mathbf A}{\mathbf D}^{-1}+{\mathbf D}{\mathbf A}^{-1}$ is

$$\medskip P(X)=\frac{0.2}{\alpha \beta \gamma \delta }(64.0+X(32.0\gamma
+40.0\delta )+X^{2}(80.0\alpha \gamma +20.0\beta \delta +20.0\gamma \delta )$$

$$ \qquad \qquad +X^{3}(50.0\alpha \gamma \delta +10.0\beta \gamma
\delta )+5.0X^{4}\alpha \beta \gamma \delta ).$$

 Since the coefficient $\frac{0.2\ast 64}{\alpha \beta \gamma
\delta }$ is positive and provides the product of all the zeroes of the characteristic
polynomial (i.e. the eigenvalues), we conclude that there is no zero eigenvalues
and consequently the considered matrix ${\mathbf A}{\mathbf D}^{-1}+{\mathbf D}{\mathbf A}^{-1}$ is nonsingular.
Therefore, applying condition (iv) of Theorem 8, we obtain that the matrix $\mathbf A$ is $D$-stable. In
order to get this result we carried out an estimated number of floating operations not greater than $50$.
It is clear that, with respect to  $k$$n!$, with $k > 0$ large, $n = 4$ (which is the number
of operations required by Hurwitz determinants methods),  the save of
operations is significant and is increasing with increasing $n$.

\bigskip
\section{Applications to dynamical systems}
\subsection{Relative stability of mechanical systems}
The linearized equations of the perturbed motion for a wide class of linear mechanical systems with $n$ degrees of freedom is as follows (see \cite{AD}, \cite{KOK}, also \cite{GANT}, \cite{KAB} for the mechanical models):
\begin{equation}\label{systmotion}
{\mathbf A}\ddot{u}(t) + {\mathbf B}\dot{u}(t) + {\mathbf C}u(t) = 0, \qquad t \geq 0,
\end{equation}
where ${\mathbf A}$, $\mathbf B$ and ${\mathbf C} \in {\mathcal M}^{n \times n}$ are mass, damping and stiffness matrices, with $\mathbf A$ being symmetric positive definite, $u(t) \in {\mathbb R}^n$ is the vector of the generalized coordinates, $\dot{u}(t)$ is the vector of the generalized speeds and $t \in {\mathbb R}$ denotes time.

Using the orthogonal similarity transformation ${\mathbf Q}$ which puts the mass matrix $\mathbf A$ into its canonical form ${\mathbf D}_{(\mathbf A)} = {\mathbf Q}^T{\mathbf A}{\mathbf Q}$, which is positive diagonal, we transform System \eqref{systmotion} as follows:
\begin{equation}\label{systmotion2}
{\mathbf D}_{({\mathbf A})}\ddot{u}(t) + \widetilde{{\mathbf B}}\dot{u}(t) + \widetilde{{\mathbf C}}u(t) = 0, \qquad t \geq 0,
\end{equation}
where ${\mathbf D}_{({\mathbf A})} \in {\mathcal M}^{n \times n}$ is the canonical form of $\mathbf A$, $\widetilde{{\mathbf B}} = {\mathbf Q}^T{\mathbf B}{\mathbf Q}$, $\widetilde{{\mathbf C}} = {\mathbf Q}^T{\mathbf C}{\mathbf Q}$.

The following definitions, which are of a certain motivation and mechanical background were provided in \cite{KOK}.

{\bf Definition 11}. System \eqref{systmotion} is called {\it $D$-stable} if the system
\begin{equation}\label{systper1}
{\mathbf A}\ddot{u}(t) + {\mathbf D}({\mathbf B}\dot{u}(t) + {\mathbf C}u(t)) = 0
\end{equation}
is asymptotically (Lyapunov) stable for every positive diagonal matrix $\mathbf D$.

{\bf Definition 12}. System \eqref{systmotion} is called {\it additive $D$-stable with respect to the coordinates} if the system
\begin{equation}\label{systper2}
{\mathbf A}\ddot{u}(t) + {\mathbf B}\dot{u}(t) + ({\mathbf C} + {\mathbf D})u(t) = 0
\end{equation}
is asymptotically stable for every nonnegative diagonal matrix $\mathbf D$.

{\bf Definition 13}. System \eqref{systmotion} is called {\it additive $D$-stable with respect to the speeds} if the system
\begin{equation}\label{systper3}
{\mathbf A}\ddot{u}(t) + ({\mathbf B}+ {\mathbf D})\dot{u}(t) + {\mathbf C}u(t) = 0
\end{equation}
is asymptotically stable for every nonnegative diagonal matrix $\mathbf D$.

{\bf Definition 14}. System \eqref{systmotion} is called {\it additive $D$-stable} if it is additive $D$-stable with respect to both the speeds and the coordinates.

Using first-order formalism (see \cite{AD}), we transform System \eqref{systmotion} into the following form:
\begin{equation}\label{systtrans}
\dot{x}(t) = \widetilde{{\mathbf A}}x(t),
\end{equation}
where $\widetilde{{\mathbf A}}$ is an $2n \times 2n$ matrix of the form
\begin{equation}
\widetilde{{\mathbf A}} = \begin{pmatrix} -{\mathbf A}^{-1}{\mathbf B} & & -{\mathbf A}^{-1}{\mathbf C} \\
{\mathbf I} & & {\mathbf O} \end{pmatrix},
\end{equation}
and $x(t) \in {\mathbb R}^{2n}$:
$$x(t) = \begin{pmatrix}\dot{u}(t) \\ u(t) \end{pmatrix}. $$

Asymptotic stability of System \eqref{systtrans} is equivalent to stability of the matrix $\widetilde{{\mathbf A}}$. Now we analyze the conditions of $D$-stability of System \eqref{systmotion} in terms of the properties of $\widetilde{{\mathbf A}}$.

It is easy to see, that $D$-stability of System \eqref{systmotion} is equivalent to stability of all the matrices $\widehat{{\mathbf D}}\widetilde{{\mathbf A}}$, where $\widehat{{\mathbf D}}$ is a block diagonal matrix of the form
$$ \widehat{{\mathbf D}} = \begin{pmatrix} {\mathbf D} & {\mathbf O} \\
{\mathbf O} & {\mathbf I} \end{pmatrix},$$
$\mathbf D$ is an $n \times n$ positive diagonal matrix. Note that the above condition is weaker than $D$-stability.

Considering the system of the form \eqref{systmotion2}, we put instead of ${\mathbf A}$ its positive diagonal canonical form ${\mathbf D}_{({\mathbf A})}$. Thus
$$ \widetilde{{\mathbf A}} = \begin{pmatrix} -{\mathbf D}_{({\mathbf A})}^{-1}\widetilde{{\mathbf B}} & & -{\mathbf D}_{({\mathbf A})}^{-1}\widetilde{{\mathbf C}} \\
{\mathbf I} & & {\mathbf O} \end{pmatrix} = \begin{pmatrix} {\mathbf D}_{({\mathbf A})}^{-1} & {\mathbf O} \\
{\mathbf O} & {\mathbf I} \end{pmatrix}\begin{pmatrix} -\widetilde{{\mathbf B}} & -\widetilde{{\mathbf C}} \\
{\mathbf I} & {\mathbf O} \end{pmatrix}.$$

It follows that $D$-stability of the matrix $\widetilde{{\mathbf A}}$ implies the asymptotic stability of System \eqref{systmotion2} for arbitrary values of the masses of its elements (see \cite{KOK}).

The additive $D$-stability of System \eqref{systmotion} with respect to the speeds and to the coordinates means that the matrices $ \begin{pmatrix} -{\mathbf A}^{-1}{\mathbf B}+ {\mathbf D} & & -{\mathbf A}^{-1}{\mathbf C}\\
{\mathbf I} & & {\mathbf O} \end{pmatrix},
$ and $\begin{pmatrix} -{\mathbf A}^{-1}{\mathbf B} & & -{\mathbf A}^{-1}{\mathbf C} + {\mathbf D}\\
{\mathbf I} & & {\mathbf O} \end{pmatrix},
$ respectively, are stable for any nonnegative diagonal matrix $\mathbf D$. The additive $D$-stability with respect to the speeds means that  the matrix
$$\widetilde{{\mathbf A}} + \begin{pmatrix} {\mathbf D} & {\mathbf O} \\
{\mathbf O} & {\mathbf O} \end{pmatrix}$$
preserves stability for all nonnegative diagonal matrices $\mathbf D$. This condition is obviously weaker than additive $D$-stability of $\widetilde{{\mathbf A}}$. Thus we conclude, that additive $D$-stability of $\widetilde{{\mathbf A}}$ implies the additive $D$-stability of System \eqref{systmotion} with respect to the speeds.

Furthermore, additive $D$-stability of System \eqref{systmotion} with respect to the coordinates means that the matrix
$$\begin{pmatrix} -{\mathbf A}^{-1}{\mathbf B} & -{\mathbf A}^{-1}{\mathbf C} + {\mathbf D}\\
{\mathbf I} & {\mathbf O} \end{pmatrix} = \widetilde{{\mathbf A}} + \begin{pmatrix} {\mathbf O} & {\mathbf D}\\
{\mathbf O} & {\mathbf O} \end{pmatrix}$$
preserves stability for all nonnegative diagonal matrices $\mathbf D$.

In its turn, additive $D$-stability of System \eqref{systmotion} means that the matrix
$$\widetilde{{\mathbf A}} + \begin{pmatrix} {\mathbf D} & {\mathbf D}\\
{\mathbf O} & {\mathbf O} \end{pmatrix}$$
preserves stability for all nonnegative diagonal matrices $\mathbf D$.

However, preserving stability under the perturbations of Forms \eqref{systper1}-\eqref{systper3}, System \eqref{systmotion} may completely lose its transient response properties. To avoid this case, we should consider other characteristics of System \eqref{systmotion}, such as {\it damping ratio} (for the definition and properties see \cite{DORB}, \cite{GANS}, \cite{ROMI}). The minimal damping ratio describes the oscillation behavior of System \eqref{systmotion}: as it decreases, the transient response of \eqref{systmotion} becomes increasingly oscillatory (see \cite{DORB}, p. 310). A variety of engineering problems requires to bound the minimal damping ratio in order to improve the transient response. Geometrically, it means that all the eigenvalues of System \eqref{systmotion} should be placed in the sector ${\mathbb C}^-_{\theta}$ with the inner angle $2\theta$ around the negative directions of the real axis. Here, the minimal damping ratio $\zeta = \cos(\theta)$ (see, e.g. \cite{GUJU}, \cite{DBE}). The value $\theta = \frac{\pi}{4}$ often arises in engineering practise (see \cite{SON}, p. 16). Robust aspects of this problem are studied in \cite{BARMI} from the polynomial point of view. Here, we analyze the cases when the minimal damping ratio (so-called relative stability) is preserved under specific perturbations. We introduce the following definition.

{\bf Definition 15}. System \eqref{systmotion} is called {\it relatively $D$-stable with the minimal damping ratio $\zeta$} if the minimal damping ratio of the perturbed system
$${\mathbf A}\ddot{u}(t) + {\mathbf D}({\mathbf B}\dot{u}(t) + {\mathbf C}u(t)) = 0$$
is less than $\zeta$ for every positive diagonal matrix $\mathbf D$.

Using the crossway from System \eqref{systmotion} to System \eqref{systtrans}, we obtain that the conditions of relative $D$-stability with the minimal damping ratio $\zeta$ are equivalent to ${\mathbb C}^-_{\theta}$-stability ($\zeta = \cos(\theta)$) of all the $2n \times 2n$ matrices $\widehat{{\mathbf D}}\widetilde{{\mathbf A}}$, where $\widehat{{\mathbf D}}$ is a block diagonal matrix of the form
$$ \widehat{{\mathbf D}} = \begin{pmatrix} {\mathbf D} & {\mathbf O} \\
{\mathbf O} & {\mathbf I} \end{pmatrix},$$
$\mathbf D$ is an $n \times n$ positive diagonal matrix. It means, that all the eigenvalues of the perturbed matrix
$$\widehat{{\mathbf D}}\widetilde{{\mathbf A}} = \begin{pmatrix} {\mathbf D} & {\mathbf O} \\
{\mathbf O} & {\mathbf I} \end{pmatrix}\begin{pmatrix} -{\mathbf A}^{-1}{\mathbf B} & & -{\mathbf A}^{-1}{\mathbf C}\\
{\mathbf I} & & {\mathbf O} \end{pmatrix}$$
should stay in ${\mathbb C}^-_{\theta}$ for any block diagonal matrix $\widehat{\mathbf D}$. This condition is obviously weaker than $({\mathbb C}^-_{\theta}, D)$-stability of $\widetilde{{\mathbf A}}$. Thus we conclude, that $({\mathbb C}^-_{\theta}, D)$-stability of $\widetilde{{\mathbf A}}$ implies the relative $D$-stability with the minimal damping ratio $\zeta$ of System \eqref{systmotion}.

\subsection{Stability of fractional-order systems}
Consider a linear system in the following form:
\begin{equation}\label{sysf}d^{\gamma}x(t) = {\mathbf A}x(t), \end{equation}
with $0 < \gamma \leq 2$, $x(0) = x_0$. It is known (see \cite{MAT}, \cite{MAT1}, \cite{MOS}, \cite{SMF}) that System \eqref{sysf} is asymptotically stable if and only if all eigenvalues $\lambda$ of ${\mathbf A}$ satisfy the inequality $|\arg(\lambda)|> \gamma\dfrac{\pi}{2}$. For $1 \leq\gamma < 2$, this corresponds to ${\mathfrak D}$-stability with respect to the stability region ${\mathfrak D} = {\mathbb C}^-_{\theta}$, where $\theta = \pi(1 - \frac{\gamma}{2})$.

Let the system matrix $\mathbf A$ be $({\mathfrak D}, {\mathcal D})$-stable with respect to ${\mathfrak D} = {\mathbb C}^-_{\theta}$. Then each system of the perturbed family
\begin{equation}\label{sysfper}d^{\gamma}x(t)  = ({\mathbf D}{\mathbf A})x(t), \end{equation}
with $\gamma = \frac{2(\pi - \theta)}{\pi}$, is asymptotically stable for every positive diagonal matrix $\mathbf D$.

\subsection{Some classical models}

Here we introduce some examples in applications where the concept of
$D$-stability plays an important role. In
particular, we point out an example related to a system of differential
equations exhibiting different time scales, presented in \cite{JOHNT} and a detailed simple example related to a general
economic problem reported in \cite{GIZ}. More, the applications of $D$-stability to the local stability analysis of steady states in the models of biochemical reaction networks (namely, to the analysis of the Glansdorff--Prigogine criterion, which gives a sufficient condition for the local stability) is considered in \cite{EGS}.

Firstly, we recall that this notion arises naturally in problems exhibiting
different time scales. \ In fact, consider a problem of the form%
\begin{equation}\label{systmod}%
\begin{array}
[c]{c}%
\varepsilon_{1}x_{1}^{\prime}=f_{1}(x_{1},...,x_{n})\\
\varepsilon_{2}x_{2}^{\prime}=f_{2}(x_{1},...,x_{n})\\
\vdots\\
\varepsilon_{n}x_{n}^{\prime}=f_{n}(x_{1},...,x_{n})
\end{array}
\end{equation}

\noindent where $f_{i}(0,..0)=0,i=1,...,n.$ Let ${\mathbf A}$ be the $n\times n$ matrix
obtained by linearizing this differential system at the origin. Then the origin is a linearly stable equilibrium
for all positive values of parameters \ $\varepsilon_{1},...,\varepsilon_{n}$
if and only if ${\mathbf A}$ is $D$-stable.\medskip

In order to cite another example, we resort to economic analysis. For
instance, the concept of $D$-stability is involved in the following economic
model. Let $x_{_{1}},...,x_{n}$ be the values of $n$ economic variables (e.g.
market prices of $n$ goods exchanged in a competitive market) and let $a$ be a
shift parameter. \ We assume that there are $n$ functional relations linking
the $x_{i}$'s and $a$, with the $i$th relation given by $f_{i}(x_{_{1}%
},...,x_{n},a).$ For a given value of $a=a^{\ast},$ an \textit{equilibrium
position} is defined as a set of values $\overline{x}_{_{1}},...,\overline
{x}_{n}$ such that
\begin{equation}
f_{i}(\overline{x}_{_{1}},...,\overline{x}_{n},a^{\ast})=0\text{ \ \ \ \ \ for
}i=1,...,n
\end{equation}

Assume that a small change in $a$ occurs. Then the change in the equilibrium
values of the variables is obtained by differentiating \eqref{systmod} totally to obtain
\begin{equation}\label{systmod2}
{\displaystyle\sum\limits_{j=1}^{n}}
\text{ }f_{ij}\frac{d\overline{x}_{j}}{da}=-f_{ia\text{ }}\text{ \ \ \ \ for
}i=1,...,n,
\end{equation}
where $f_{ij}=\frac{\partial f_{i}}{\partial x_{j}}$ and $f_{ia}%
=\frac{\partial f_{i}}{\partial a}$ are the partial derivatives evaluated at
the equilibrium position $(\overline{x}_{_{1}}, \ ..., \ \overline{x}_{n},a^{\ast
}).$ In matrix notation, \eqref{systmod2} can be rewritten%
$$
{\mathbf A}\left[  \frac{d\overline{x}_{j}}{da}\right]  =-\left[  f_{ia\text{ }}\right],
$$
where the bracketed terms are $n\times1$ vectors and $\mathbf A$ is a square
$n\times n$ matrix with $(i,j)-$element equal to $\ f_{ij}.$

In addition it is a assumed that a dynamic adjustment process operates to
determine the time paths of variables $x_{_{1}},...,x_{n}$ \ when the system
is out of equilibrium. We write this adjustment process as follows:%
\begin{equation}\label{systmod3}
\frac{dx_{i}}{dt}=b_{i}\text{ }f_{i}(x_{_{1}},...,x_{n},a^{\ast})\text{
\ \ \ \ for }i=1,...,n,
\end{equation}
where $b_{i}$ are positive constants; each of them is called the "speed of
adjustment" for the $i$th variable. In particular, in a sufficiently small
neighborhood of equilibrium, the adjustment process can be approximated by the
linear terms of a Taylor expansion, so that \eqref{systmod3} becomes
\begin{equation}\label{systmod4}
\frac{dx_{i}}{dt}=b_{i}\sum_{j=1}^{n}
f_{ij}(x_{j}-\overline{x}_{j}) \quad \mbox{for} \quad i=1,...,n.
\end{equation}

Let $z_{i}=x_{i}-\overline{x}_{i}.$ \ Then \eqref{systmod4} in matrix notation can be
rewritten as%
\begin{equation}\label{systmod5}
\dot{z}= {\mathbf B}{\mathbf A}z,
\end{equation}
where $\dot{z}=\left[  \frac{dz_{i}}{dt}\right]$ and $\mathbf B$ is a
diagonal matrix with $b_{i}>0$ in the $i$th diagonal position.

Then system \eqref{systmod5} is stable if and only if the real parts of all the eigenvalues
of ${\mathbf B}{\mathbf A}$ are negative, i.e. if and only if $A$ is $D$-stable.\medskip

\section{Concluding remarks}
The concept of $D$-stability enjoyed a large inhomogeneous literature in the
past decades. Here, the challenging aim of this paper is to unify the basic
different theoretical approaches by means of a generalization of the $D$-stability
concept to a family of unbounded regions defined in terms of Linear Matrix Inequalities
(LMI regions). The novelty of this work is properly the introduction
of this promising general point of view, which allows to extend properties and
to simplify computations required for the solution of the $D$-stability characterization
problem. Actually, the problem of computational complexity is studied
in some details. Many theorems are proved which show the power of this new
approach. Moreover, together with multiplicative $D$-stability, we consider additive
$D$-stability and we show how it can be used efficiently in the study of
stability of dynamical systems.
This paper is intended to be just the first one on this subject, which in our
opinion is worthy of further studies, since it appears promising and fruitful.



%
%



\end{document}